\documentclass{amsart}

\newtheorem{theorem}{Theorem}[section]

\newtheorem{cor}[theorem]{Corollary}
\newtheorem{conj}[theorem]{Conjecture}
\newtheorem{prop}[theorem]{Proposition}

\theoremstyle{definition}

\theoremstyle{remark}
\newtheorem{remark}[theorem]{Remark}

\numberwithin{equation}{section}

\newcommand{\abs}[1]{\lvert#1\rvert}

\DeclareSymbolFont{AMSb}{U}{msb}{m}{n}
\DeclareMathSymbol{\Z}{\mathalpha}{AMSb}{"5A}

\DeclareMathSymbol{\nmid}{\mathrel}{AMSb}{"2D}

\begin{document}
\newcommand{\beqs}{\begin{equation*}}
\newcommand{\eeqs}{\end{equation*}}
\newcommand{\beq}{\begin{equation}}
\newcommand{\eeq}{\end{equation}}
\newcommand\mylabel[1]{\label{#1}}
\newcommand\eqn[1]{(\ref{eq:#1})}
\newcommand\gausspoly[2]{\begin{bmatrix} #1 \\ #2\end{bmatrix}_q}
\newcommand\mytwid[1]{\overset {\text{\lower 3pt\hbox{$\sim$}}}#1}

\title[K.~Saito's Conjecture for Nonnegative Eta Products]{K.~Saito's 
Conjecture for Nonnegative \\
Eta Products and Analogous Results\\
for Other Infinite Products}

\author{Alexander Berkovich}
\address{Department of Mathematics, University of Florida, Gainesville,
Florida 32611-8105}
\email{alexb@math.ufl.edu}          
\thanks{Both authors were supported in part by NSA Grant MSPF-06G-150.}

\author{Frank G. Garvan}
\address{Department of Mathematics, University of Florida, Gainesville,
Florida 32611-8105}
\email{frank@math.ufl.edu}          


\date{January 24, 2007}  



\begin{abstract}
We prove that the Fourier coefficients of a certain general
eta product considered by K.~Saito
are nonnegative. The proof is elementary and depends on a
multidimensional theta function identity. The $z=1$ case 
is an identity for the generating function for $p$-cores due to
Klyachko \cite{K1} and Garvan, Kim and Stanton \cite{GKS}.
A number of other infinite products are shown to have
nonnegative coefficients. In the process a new generalization
of the quintuple product identity is derived.
\end{abstract}

\maketitle

\section{Introduction} \label{sec:intro}

Throughout this paper $q=\exp(2\pi i\tau)$ with $\Im \tau >0$ so that $\abs{q}<1$. 
As usual the
Dedekind eta function is defined as
\beq
\eta(\tau) := \exp(\pi i\tau/12) \prod_{n=1}^\infty(1- \exp(2\pi in\tau))
= q^{1/24} \prod_{n=1}^\infty(1-q^n).
\mylabel{eq:etadef} 
\eeq
An eta product is a finite product of the form
\beq
\prod_k \eta(k\tau)^{e(k)},
\mylabel{eq:etaproddef} 
\eeq
where the $e(k)$ are integers. K.~Saito \cite{S1} considered eta products
that are connected with elliptic root systems and considered the
problem of determining when all the Fourier coefficients of such eta products
are nonnegative. Subsequent work contains the following 
\begin{conj} 
\label{conj1}
{\rm (K.~Saito \cite{S3})} Let $N$ be a positive integer.
The eta product
\beq
S_N(\tau):=\frac{ \eta(N\tau)^{\phi(N)}}
{\prod_{d\mid N} \eta(d\tau)^{\mu(d)}}
\mylabel{eq:setaproddef} 
\eeq
has nonnegative Fourier coefficients. 
\end{conj}
The conjecture has been proved for $N=2,3,4,5,6,7,10$ by K.~Saito \cite{S1}, \cite{S2},
\cite{S3}, \cite{S4}, \cite{S5}, for prime powers $N=p^\alpha$ by 
T.~Ibukiyama \cite{I1}, 
and for $\gcd(N,6)>1$ by K.~Saito and S.~Yasuda \cite{SY},
who also showed that for general $N$, the coefficient of $q^n$
in $S_N(\tau)$ is nonnegative for sufficiently large $n$.
We prove the conjecture for general $N$. 

The case $N=p$ (prime) occurs in the study of $p$-cores.
A partition is a $p$-core if it has no hooks of length $p$ \cite{GKS},
\cite{JK}. 
$p$-cores are important in the study of $p$-modular
representations of the symmetric group $S_n$. Define
\beq
E(q) := \prod_{n=1}^\infty (1 - q^n),
\mylabel{eq:Edef} 
\eeq
and let $a_t(n)$ denote the number of partitions of $n$ that are
$t$-cores. It is well known that for any positive integer $t$
\beq
\sum_{n\ge0} a_t(n) q^n =  \frac{E(q^t)^t}{E(q)}.
\mylabel{eq:pcore1} 
\eeq
This result is originally due to Littlewood \cite{L}.
See \cite{GKS} for a combinatorial proof. 
Thus \eqn{pcore1} implies that Conjecture \ref{conj1} holds for
$N=p$ prime since
\beq
S_p(\tau) = \frac{\eta(p\tau)^p}{\eta(\tau)}
= q^{(p^2-1)/24} \frac{E(q^p)^p}{E(q)}.
\mylabel{eq:Sp} 
\eeq
Granville and Ono \cite{GO} have proved that $a_t(n) > 0$
for all $t\ge4$ and all $n$.
We also need the following identity due to Klyachko \cite{K1}
\beq
\sum_{\substack{\vec n\in \mathbb Z^t \\ \vec n\cdot\vec{1}_t=0}} q^{\frac{t}{2}
\vec n\cdot\vec{n}+\vec b_t\cdot\vec{n}}=
\frac{E(q^t)^t}{E(q)}, 
\mylabel{eq:kid} 
\eeq
where $\vec{1}_t=(1,1,\ldots,1)\in \mathbb Z^t$, $\vec b_t=(0,1,2,\ldots,t-1)$,
and $t$ is any positive integer.
See \cite{GKS} for a combinatorial proof. See also \cite[Prop.1.29]{G94}
and \cite[\S2]{BG06b}.
Our proof of K.~Saito's Conjecture
depends on the following $z$-analogue of \eqn{kid}.
\begin{theorem}
\mylabel{thm1} Let $a\ge 2$ be an integer. Then for $z\ne0$ and $\abs{q}<1$
we have
\begin{align}
C_a(z;q)&:= 
\sum_{\substack{\vec n=(n_0,n_1,\dots, n_{a-1})\in \mathbb Z^a \\ \vec n\cdot\vec{1}_a=0}} 
q^{\frac{a}{2}\vec n\cdot\vec{n}+\vec b_a\cdot\vec{n}}
\sum_{j=0}^{a-1} z^{an_j + j}
\mylabel{eq:Cazq}\\ 
&= E(q) E(q^a)^{a-2}
\prod_{n=1}^\infty \frac{(1-z^a q^{a(n-1)})(1-z^{-a}q^{an})}
                        {(1-z q^{n-1})(1-z^{-1}q^n)}
\nonumber
\end{align}
\end{theorem}
We note that \eqn{kid} follows from \eqn{Cazq} by letting $a=t$ and $z\to1$.
The case $a=3$ is equivalent to \cite[(1.23)]{HGB}. See \cite[\S3.3]{BG06a}.
The case $a=2$ can be written as
\beq
\sum_{n\in\mathbb Z} q^{2n^2+n}\left(z^{2n} + z^{2n+1}\right)
= 
\prod_{n=1}^\infty (1+z q^{(n-1)})(1+z^{-1}q^{n})(1 - q^n)
\mylabel{eq:Cazq2}, 
\eeq
which follows easily from Jacobi's triple product identity
\cite[(2.2.10)]{Andbook}.
Klyachko \cite{K1} proved \eqn{kid} by showing that it was a
special case of the Macdonald identity \cite{M1} for the affine root
system $A_{t-1}$. Hjalmar Rosengren \cite{R1} has observed that
Theorem \ref{thm1} (i.e. \eqn{Cazq}) is also a special
case of the Macdonald identity for $A_{a-1}$.

If we rewrite the right side of \eqn{kid} in terms of eta products,
and apply Jacobi's transformation $\tau \mapsto -1/\tau$ then
we are led to the identity
\beq
\sum_{\substack{\vec n\in \mathbb Z^t \\ \vec n\cdot\vec{1}_t=0}} 
\omega_t^{\vec b_t\cdot\vec{n}} q^{\frac{1}{2} \vec n\cdot\vec{n}}
=
\frac{E(q)^t}{E(q^t)},
\mylabel{eq:ckid} 
\eeq
where $\omega_t:=\exp(2\pi i/t)$. The proof uses well known transformation
formulas for the eta function \cite[Thm.3.1]{Ap} and multidimensional
theta functions \cite[(5),p.205]{Sch}.
See \cite[Prop.2.29]{G94} for an elementary proof. Equation \eqn{ckid}
has the following $z$-analogue. 
\begin{theorem}
\mylabel{thm2} Let $a$ and $j$ be integers where
$a\ge2$ and $0\le j \le a-1$. Then for $z\ne0$ and $\abs{q}<1$
we have
\begin{align}
B_{j,a}(z;q)&:=
\sum_{\substack{\vec n=(n_0,n_1,\dots, n_{a-1})\in \mathbb Z^a \\ 
\vec n\cdot\vec{1}_a=0}}
z^{n_j} 
\omega_a^{\vec b_a\cdot\vec{n}}
q^{\frac{1}{2}\vec n\cdot\vec{n}}
\mylabel{eq:Bjazq}\\ 
&= E(q)^{a-2} E(q^a)
\prod_{n=1}^\infty \frac{(1-z q^{n-1})(1-z^{-1}q^{n})}
                        {(1-z q^{a(n-1)})(1-z^{-1}q^{an})}
\nonumber
\end{align}
\end{theorem}
Clearly, \eqn{ckid} follows from \eqn{Bjazq} by letting $a=t$ and $z\to 1$.
Again the $a=2$ case follows easily from Jacobi's triple
product identity.

\subsection*{Notation} 
We use the following notation for finite products
$$
(z;q)_n=(z)_n=
\begin{cases}
\prod_{j=0}^{n-1}(1-zq^j), & n>0 \\
1,                         & n=0.
\end{cases}
$$
For infinite products we use
$$
(z;q)_\infty=(z)_\infty = \lim_{n\to\infty} (z;q)_n
=\prod_{n=1}^\infty (1-z q^{(n-1)}),
$$
$$
(z_1,z_2,\dots,z_k;q)_\infty = (z_1;q)_\infty (z_2;q)_\infty
\cdots (z_k;q)_\infty,
$$
$$
[z;q]_\infty = (z;q)_\infty (z^{-1}q;q)_\infty=
\prod_{n=1}^\infty (1-z q^{(n-1)})(1-z^{-1}q^{n}),
$$
$$
[z_1,z_2,\dots,z_k;q]_\infty = [z_1;q]_\infty [z_2;q]_\infty
\cdots [z_k;q]_\infty,
$$
for $\abs{q}<1$ and $z$, $z_1$, $z_2$,\dots, $z_k\ne 0$.


\section{Proof of Theorems \ref{thm1} and \ref{thm2}}
\label{sec:proofthm1}
Suppose $a\ge2$. 
The idea is to show both sides of \eqn{Cazq} satisfy the
same functional equation as $z\to zq$ and agree for enough
values of $z$.
Define
\beq
R_a(z;q) =  E(q) E(q^a)^{a-2}
            \frac{[z^a;q^a]_\infty}{[z;q]_\infty},
\mylabel{eq:Radef} 
\eeq
which is the right side of \eqn{Cazq}.
An easy calculation gives
\beq
R_a(zq;q) = z^{-(a-1)} R_a(z;q).
\mylabel{eq:Rafe} 
\eeq
We show that basically the $a$ terms in the definition of $C_a(z;q)$
are permuted cyclically as $z\to zq$. To this end we define
\beq
Q_a(\vec{n}) = \frac{a}{2}\vec n\cdot\vec{n}+\vec b_a\cdot\vec{n},
\mylabel{eq:Qadef} 
\eeq
\beq
F_j(z;q) :=
\sum_{\substack{\vec n=(n_0,n_1,\dots, n_{a-1})\in \mathbb Z^a \\ \vec n\cdot\vec{1}_a=0}} 
 z^{an_j + j} q^{Q_a(\vec{n})}
\qquad (0\le j \le a-1).
\mylabel{eq:Fjdef} 
\eeq

Now suppose $1 \le j \le a-1$.
Let $\vec{e}_0=(1,0,\dots,0)$, 
$\vec{e}_1=(0,1,\dots,0)$, \dots,
$\vec{e}_{a-1}=(0,0,\dots,0,1)$ be the standard unit vectors,
$\vec{n}=(n_0,n_1,\dots, n_{a-1})\in \mathbb Z^a$,
and
$\vec{n}\,'=(n_1,n_2,\dots, n_{a-1},n_0) + \vec{e}_{j-1} - \vec{e}_{a-1}$.
An easy calculation gives
\beq
Q_a(\vec{n}\,') - Q_a(\vec{n})
= a n_j + j - \vec{n}\cdot\vec{1}_a.
\mylabel{eq:Qdiff} 
\eeq
Hence
\begin{align}
F_{j-1}(z;q) &= \sum_{\substack{\vec n\in \mathbb Z^a \\ 
\vec n\cdot\vec{1}_a=0}}
z^{an_{j-1}+(j-1)} q^{Q_a(\vec{n})}\mylabel{eq:Fjm1}\\ 
&= \sum_{\substack{\vec{n}\,'\in \mathbb Z^a \\ 
\vec n\,'\cdot\vec{1}_a=0}}
z^{a(n_{j}+1)+(j-1)} q^{Q_a(\vec{n}\,')}\nonumber\\
&= \sum_{\substack{\vec{n}\in \mathbb Z^a \\ 
\vec n\cdot\vec{1}_a=0}}
z^{an_{j}+j+(a-1)} q^{Q_a(\vec{n})+an_j+j}\nonumber\\
&=z^{(a-1)} F_j(zq;q),
\nonumber
\end{align}
and
\beq
F_j(zq;q) = z^{-(a-1)} F_{j-1}(z;q).
\mylabel{eq:Fjfe} 
\eeq
Similarly we find that
\beq
F_0(zq;q) = z^{-(a-1)} F_{a-1}(z;q),
\mylabel{eq:F0fe} 
\eeq
by using the result that
\beq
Q_a(\vec{n}\,') - Q_a(\vec{n})
= a n_{0} - \vec{n}\cdot\vec{1}_a,
\mylabel{eq:Qdiff0} 
\eeq
where
$\vec{n}\,'=(n_1,n_2,\dots, n_{a-1},n_0)$.

Since
\beq
C_a(z;q) = \sum_{j=0}^{a-1} F_j(z;q),
\mylabel{eq:Casum} 
\eeq
we have
\beq
C_a(zq;q) = z^{-(a-1)} C_a(z;q).
\mylabel{eq:Cafe} 
\eeq
In view of \cite[Lemma 2]{ASD} or \cite[Lemma 1]{HGB}
it suffices to show that \eqn{Cazq} holds for $a$ distinct
values of $z$ with $\abs{q} < \abs{z} \le 1$.
It is clear that
\beq
C_a(z;q) = R_a(z;q) = 0,
\mylabel{eq:Cazqzero} 
\eeq
for $z=\exp(2\pi ik/a)$ for $1 \le k \le a-1$.
Finally, \eqn{Cazq} holds for $z=1$ since
\beq
C_a(1;q) = a 
\sum_{\substack{\vec n\in \mathbb Z^a \\ \vec n\cdot\vec{1}_a=0}} 
q^{\frac{a}{2}
\vec n\cdot\vec{n}+\vec b_a\cdot\vec{n}}
= a\frac{E(q^a)^a}{E(q)}=R_a(1;q),
\eeq
by \eqn{kid} with $t=a$.
This completes the proof of Theorem \ref{thm1}.

The proof of Theorem \ref{thm2} is similar.
First by considering a cyclic permutation one can show that the
definition of $B_{j,a}(z;q)$ is
independent of $j$. Next one show that both sides of
\eqn{Bjazq} satisfy the same functional equation
\beq
\Phi_{a}(zq^a;q) = q^{-\binom{a}{2}} (-z)^{-(a-1)} \Phi_{a}(z;q).
\mylabel{eq:Phitrans}
\eeq
For the right side this is easy. For left side one uses the
transformation
$\vec{n}\mapsto \vec{n}\,'=\vec{n} + \vec{1}_a - a \vec{e}_j$.
Since both sides satisfy the same functional equation \eqn{Phitrans}
and are analytic for $z\ne0$, one needs only to verify that
the identity holds for $a$ distinct values of $z$ on the region
$\abs{q}^a < \abs{z} \le 1$. For $1\le k \le a-1$, one
uses the transformation 
$\vec{n}\mapsto \vec{n}\,'=\vec{n} + k (\vec{e}_{j+1}-\vec{e}_j)$
to find that
\beq
B_{j,a}(q^k;q)= \omega_a^k B_{j+1,a}(q^k;q)
= \omega_a^k B_{j,a}(q^k;q),
\mylabel{eq:Bzero0}
\eeq
so that
\beq
B_{j,a}(q^k;q)=0.
\mylabel{eq:Bzero1}
\eeq
Hence, both sides of \eqn{Bjazq} are zero for $z=q^k$ ($1\le k \le a$)
and agree at $z=1$ by \eqn{ckid} with $t=a$. The theorem follows.

\section{Proof of K.~Saito's Conjecture} 
\label{sec:proofSaito}

First we show that
\beq
\prod_{d\mid M} E(q^d)^{\mu(d)} =
\prod_{\substack{ n\ge 1\\ (n,M)=1}} (1 - q^n),
\mylabel{eq:Eprop} 
\eeq
for any positive integer $M$.
Now
\beq
\prod_{d\mid M} E(q^d)^{\mu(d)}
=
\prod_{d\mid M} \prod_{m=1}^\infty (1- q^{dm})^{\mu(d)}
=
\prod_{n=1}^\infty (1- q^{n})^{\varepsilon(n)},
\mylabel{eq:Eep} 
\eeq
where
\beq
\varepsilon(n) =
\sum_{d\mid M\,\&\,d\mid n} \mu(d)
=
\sum_{d\mid (M,n)} \mu(d)
=
\begin{cases}
1 & \mbox{if $(M,n)=1$}\\
0 & \mbox{otherwise},
\end{cases}
\mylabel{eq:epsimp} 
\eeq
by a well known property of the M\"obius function, and we have
\eqn{Eprop}.

For any positive integer $N$ we define
\beq
\mytwid{S}_N(q) := \frac
{E(q^N)^{\phi(N)}}
{\prod_{d\mid N} E(q^d)^{\mu(d)}}.
\mylabel{eq:SNdef} 
\eeq
We wish to show that all coefficients in the $q$-expansion of $\mytwid{S}_N(q)$
are nonnegative. We consider three cases.

\noindent
{\it Case 1.} $N=p^\alpha$ where $p$ is prime. This case was proved
by Ibukiyama \cite{I1}. Alternatively, the case $\alpha=1$ follows from 
\eqn{pcore1} and then use an easy induction on $\alpha$.

\noindent
{\it Case 2.} $N=p\,M$, where $p$ is prime, $M$ is odd and $p\nmid M$.
We have
\beq
\prod_{d\mid N} E(q^d)^{\mu(d)} =
\prod_{d\mid M} E(q^d)^{\mu(d)}  E(q^{pd})^{\mu(pd)} =
\prod_{d\mid M} \left(
\frac{E(q^d)}{E(q^{pd})}\right)^{\mu(d)}.
\mylabel{eq:Eprod} 
\eeq
By \eqn{Eprop} we have
\beq
\prod_{d\mid M} E(q^d)^{\mu(d)} =
\prod_{\substack{ n\ge 1\\ (n,M)=1}} (1 - q^n) =
\prod_{n\ge0}
\prod_{\substack{(r,M)=1\\1\le r \le M-1}} (1 - q^{Mn+r}) =
\prod_{\substack{(r,M)=1\\1\le r \le \frac{M-1}{2}}} [q^r;q^M]_\infty.
\mylabel{eq:Eprod2} 
\eeq
Now for $a$ a positive integer, $\abs{q}<1$ and $z\ne0$ we let
\beq
D_a(z;q) := \frac{E(q^a)^a}{E(q)} C_a(z;q)
         = E(q^a)^{2a-2} \frac{[z^a;q^a]_\infty}{[z;q]_\infty}
\mylabel{eq:Dazq}, 
\eeq
so that
\begin{align}
\prod_{\substack{(r,M)=1\\1\le r \le \frac{M-1}{2}}} D_p(q^r;q^M) 
&= \left(E(q^{pM})^{2p-2}\right)^{\phi(M)/2}
\prod_{\substack{(r,M)=1\\1\le r \le \frac{M-1}{2}}} 
\frac{[q^{pr};q^{pM}]_\infty}{[q^{r};q^{M}]_\infty} \nonumber\\
&=E(q^N)^{\phi(N)}
\prod_{d\mid M} \frac{ E(q^{pd})^{\mu(d)}}{ E(q^{d})^{\mu(d)}}
\qquad\mbox{(by \eqn{Eprod2})}
\mylabel{eq:Dprod}\\ 
&= \frac
{E(q^N)^{\phi(N)}}
{\prod_{d\mid N} E(q^d)^{\mu(d)}} 
\qquad\mbox{(by \eqn{Eprod})}
\nonumber\\
&= \mytwid{S}_N(q).
\nonumber
\end{align}
K.~Saito's Conjecture holds in this case since each $C_p(q^r;q^M)$
has nonnegative coefficients by Theorem \ref{thm1}, and
$E(q^{pM})^p/E(q^M)$ has nonnegative coefficients by \eqn{pcore1}
so that each $D_p(q^r;q^M)$ has nonnegative coefficients.

\noindent
{\it Case 3.} $N=p^\alpha\,M$, where $p$ is prime, $M$ is odd, $p\nmid M$,
and $\alpha\ge2$. We let $N'=p\,M$. It is clear that
\beq
\prod_{d\mid N} E(q^d)^{\mu(d)}= \prod_{d\mid N'} E(q^d)^{\mu(d)}.
\mylabel{eq:Eprop2} 
\eeq
Hence                              
\beq
\mytwid{S}_N(q) = \frac{E(q^N)^{\phi(N)}}
{E(q^{N'})^{\phi(N')}} \mytwid{S}_{N'}(q) 
= \left(\frac{E(q^{p^{\alpha-1}N'})^{p^{\alpha-1}}}{E(q^{N'})}
\right)^{ (p-1)\phi(M)} \mytwid{S}_{N'}(q).
\mylabel{eq:Sprop} 
\eeq
Here $\mytwid{S}_N(q)$ is the product of two terms. The second term
$\mytwid{S}_{N'}(q)$ has 
nonnegative coefficients from Case 2. The first term has
nonnegative coefficients
using \eqn{pcore1} with $q$ replaced
with $q^{N'}$ and $t=p^{\alpha-1}$. Thus K.~Saito's Conjecture holds
in this case.

\section{Other Products with Nonnegative Coefficients} 
\label{sec:otherprods}

In this section we prove a number of results for coefficients
of other infinite products. 
For a formal power series
$$
F(q) := \sum_{n=0}^\infty a_n q^n \in \Z[[q]]
$$
we write
$$
F(q) \succeq 0,
$$
if
$a_n \ge 0$ for all $n\ge0$.
For a formal power series $F(z_1,z_2,\dots,z_n;q)$ in more than one variable
we interpret $F(z_1,z_2,\dots,z_n;q)\succeq0$ in the natural way.
The following result follows from the $q$-binomial theorem
\cite[Thm2.1]{Andbook}.
\begin{prop}
\label{atineq}If $\abs{q}$, $\abs{t} < 1$ then
\beq
\frac{ (at;q)_\infty}{(a;q)_\infty (t;q)_\infty}
=\sum_{n=0}^\infty 
\frac{t^n}{(aq^n;q)_\infty (q)_n} \succeq 0.
\mylabel{eq:atq} 
\eeq
\end{prop}

\begin{cor}
\label{atcor}

If $a$, $b$, $M$ are positive integers then
\beq
\prod_{n=0}^\infty
\frac{ (1 - q^{Mn+a+b}) }{ (1-q^{Mn+a}) (1-q^{Mn+b}) } \succeq 0.
\mylabel{eq:coratq1} 
\eeq
\end{cor}

Proposition \ref{atineq} has a finite analogue. 
For $0 \le m \le n$ the Gaussian polynomial 
\cite[p.33]{Andbook} is defined by
\beq
\gausspoly{n+m}{m} =
\frac{ (q)_{m+n}}{(q)_n (q)_m}=
\frac{ (1-q^{n+1}) \cdots (1-q^{n+m})}{(q)_m}
\mylabel{eq:gpdef} 
\eeq
Since it is 
the generating function for partitions with at most $m$ parts each $\le n$
it is a polynomial (in $q$) with positive integer coefficients.
We have
\begin{prop}
\label{atfinite}
If $L\ge0$ then
\beq
\frac{ (z_1 z_2;q)_L }{ (z_1;q)_L (z_2;q)_L }
=
\sum_{j=0}^L
\gausspoly{L}{j} \frac{ z_1^j }
{(z_1q^{L-j};q)_j (z_2q^{j};q)_{L-j}} \succeq 0.
\mylabel{eq:atqfin} 
\eeq
\end{prop}
This proposition follows from \cite[Ex1.3(i), p.20]{GR}.

The case $(t,a)=(z,qz^{-1})$ of Proposition \ref{atineq} is 
\beq
\frac{ (q;q)_\infty}
     {(z;q)_\infty (q/z;q)_\infty}=
\frac{E(q)}
     {[z;q]_\infty}
=
\sum_{n=0}^\infty
\frac{z^n}{(q)_n (z^{-1}q^{n+1};q)_\infty} \succeq 0,
\mylabel{eq:zq}             
\eeq
and is 
related to the
crank of partitions \cite{AG}. See also \cite[Eq.(5.7),p.43]{G88}.
The crank of a partition is the largest part if the partition
contains no ones, and is otherwise
the number of parts larger than the number of ones minus the number of ones.
Let $M(m,n)$ denote the number of partitions of $n$ with crank $m$.
Then
\begin{align}
(1-z) \frac{E(q)}{[z;q]_\infty}
&= \prod_{n=1}^\infty \frac{(1 - q^n)}{(1-zq^n)(1-z^{-1}q^n)} 
\mylabel{eq:crankgen} 
\\
&= 1 + (z + z^{-1} - 1)q + 
\sum_{n\ge2}\left(\sum_{m=-n}^n M(m,n) z^m\right) q^n.
\nonumber                       
\end{align}
This result follows from (1.11) and Theorem 1 in \cite{AG}.
We note the coefficients on the right side of \eqn{crankgen} are
nonnegative except for the coefficient of $z^0q^1$.
By observing that
$$
(1+z+z^2 + \cdots + z^{m-1})(z + z^{-1} - 1) = z^{-1}
+ \sum_{j=1}^{m-2} z^j + z^m \qquad(m\ge2)
$$
we have
\begin{prop}
\label{aci}
If $\abs{q}<1$, $z\ne0$ and $m\ge2$ then
\beq
(1-z^m) \frac{E(q)}{[z;q]_\infty}
= (1+z+z^2 + \dots + z^{m-1}) \prod_{n=1}^\infty \frac{(1 - q^n)}{(1-zq^n)(1-z^{-1}q^n)}
\succeq 0.
\mylabel{eq:aci} 
\eeq
\end{prop}

We will also need
\beq
\frac{E(q^t)^t}{E(q)} \succeq 0,\quad\mbox{for any positive integer $t$.}
\mylabel{eq:tcore}
\eeq
This follows from\eqn{pcore1}.

The quintuple product identity \cite[Ex5.6, p.134]{GR} can
be written as
\beq
\frac{[z^2;q]_\infty E(q)}{[z,z^3;q]_\infty}
=  \frac{E(q^3)}{[z^3,q^2z^3;q^3]_\infty} +
 z\,\frac{E(q^3)}{[z^3,qz^3;q^3]_\infty}.
\mylabel{eq:quin}
\eeq
Ekin \cite{E1} used this form of the
quintuple product identity to prove a number of inequalities for
the crank of partitions mod $7$ and $11$.
In the following Proposition we give a generalization of the
quintuple product identity. The case $a=2$ is \eqn{quin}.
\begin{theorem}
\label{GQPI} 
{\rm (Generalization of the Quintuple Product Identity)}
Suppose $a$ is a positive integer, $\abs{q}<1$ and $z\ne0$.
Then
\beq
\frac{ [z^a;q]_\infty }
     { [z,z^{a+1};q]_\infty}
=
\frac{E(q^{a+1})^2}
     {E(q)^2}
\sum_{j=0}^{a-1}
z^j \, \frac{ [q^{a-j};q^{a+1}]_\infty}
            {[z^{a+1},z^{a+1}q^{a-j};q^{a+1}]_\infty}.
\mylabel{eq:gqpi}
\eeq
\end{theorem}
\begin{proof}
We rewrite \eqn{gqpi} as
\begin{align}
&\frac{ E(q)^2 }
     {E(q^{a+1})^2}
\frac{ [z^a;q]_\infty }
     { [z;q]_\infty} \mylabel{eq:gqpib}\\
&\quad=
[z^{a+1}q,z^{a+1}q^2, \dots, z^{a+1}q^a;q^{a+1}]_\infty
\sum_{j=0}^{a-1}
z^j \, \frac{ [q^{a-j};q^{a+1}]_\infty}
            {[z^{a+1}q^{a-j};q^{a+1}]_\infty}.
\nonumber
\end{align}
Using the fact that
\beq
[zq^k;q]_\infty = (-1)^k z^{-k} q^{-\binom{k}{2}} [z;q]_\infty
\mylabel{eq:jactrans}
\eeq
it is straighforward to show that both sides of \eqn{gqpib}
satisfy the functional equation
\beq
\Phi_{a}(zq;q) = (-1)^{a-1}q^{-\binom{a}{2}} z^{1-a^2} \Phi_{a}(z;q).
\mylabel{eq:Phitrans2}
\eeq
Therefore, since both sides of \eqn{gqpib} are analytic for $z\ne0$ we
need only to verify \eqn{gqpib} for $a^2$ distinct values of $z$ in the
region $\abs{q} < z \le 1$.
For $1\le k \le a$, and $0\le n \le a$ we let
$$
z= q^{k/(a+1)} e^{2\pi i n/(a+1)},
\qquad \mbox{so that} \quad z^{a+1} = q^k \quad \mbox{and}
\quad z^a = z^{-1} q^k.
$$
We find that each term in the sum on the right side of \eqn{gqpib}
is zero except the term corresponding to $j=k-1$, and that
both sides simplify to
$$
(-1)^{k+1} q^{-\binom{k}{2}} z^{k-1} \frac{ E(q)^2}
                                         { E(q^{a+1})^2}.
$$
Thus both sides of \eqn{gqpib} agree for $a^2+a$ distinct
values of $z$ and the result follows.
\end{proof}

\begin{remark}
This theorem can also be proved using Ramanujan's
${}_1\Psi_1$-summation.
\end{remark}

\begin{remark}
Ekin's identity \cite[(38), p.287]{E1}
\beq
\frac{E(q)}{[z;q]_\infty} 
= \frac{1}{[z^2;q^2]_\infty} \sum_{n=-\infty}^\infty z^n q^{n(n-1)/2},
\mylabel{eq:EkinId1}
\eeq
implies
\beq
\frac{E(q)E(q^2)}{[z;q]_\infty}
= \frac{1}{[z^4;q^4]_\infty} \sum_{n_1,n_2=-\infty}^\infty z^{n_1+2n_2} 
q^{n_1(n_1-1)/2 + n_2(n_2-1)}\succeq	0.
\mylabel{eq:EkinId2}
\eeq
We can iterate \eqn{EkinId1} to obtain
\beq
\frac{E(q)E(q^2)E(q^4)E(q^8)\cdots}{[z;q]_\infty} \succeq 0.
\mylabel{eq:EkinIt}
\eeq
\end{remark}

\begin{cor}
\label{cor:gqpi}
Let $\abs{q}<1$ and $z\ne0$.
If $a$ and $k$ are integers with $a$ and $k\ge2$, then
\beq
(1-z^{k(a+1)})
E(q) E(q^{a+1})^{\lfloor (a+1)/2\rfloor}
\frac{[z^a;q]_\infty}
     {[z,z^{a+1};q]_\infty} \succeq 0.
\mylabel{eq:corgqpi1}
\eeq
\end{cor}
\begin{proof}

From \eqn{gqpi} we have
\begin{align}
&(1-z^{k(a+1)}) E(q) E(q^{a+1})^{\lfloor (a+1)/2\rfloor}
\frac{[z^a;q]_\infty}
     {[z,z^{a+1};q]_\infty} \nonumber \\
&=
\sum_{i=0}^{a-1}
z^i \, 
\frac{E(q^{a+1})^{\lfloor (a+1)/2\rfloor} [q^{a-i};q^{a+1}]_\infty}
            {E(q)} \mylabel{eq:gqpic} \\
&\qquad\qquad   \cdot
	(1-z^{k(a+1)})\frac{E(q^{a+1})^2}
                           {[z^{a+1},z^{a+1}q^{a-i};q^{a+1}]_\infty}.
\nonumber
\end{align}
Suppose $0 \le i \le a-1$.
By \eqn{aci} each term
\beq
(1-z^{k(a+1)})
\frac{E(q^{a+1})^2}
     {[z^{a+1},z^{a+1}q^{a-i};q^{a+1}]_\infty}\succeq 0,
\eeq
since $k\ge2$.
It remains to show that each term
\beq
\frac{E(q^{a+1})^{\lfloor (a+1)/2\rfloor} [q^{a-i};q^{a+1}]_\infty}
            {E(q)} \succeq 0.
\eeq
If $a\equiv 0\pmod{2}$ then we find that
\beq
\frac{E(q^{a+1})^{\lfloor (a+1)/2\rfloor} [q^{a-i};q^{a+1}]_\infty}
            {E(q)} 
= \prod_{\substack{ j=1\\ j\not\in\{a-i,i+1\}}}^{a/2}
\frac{E(q^{a+1})}{[q^j;q^{a+1}]_\infty} \succeq 0,
\eeq
by \eqn{zq}. If $a\equiv 1\pmod{2}$ and $i\ne (a-1)/2$ we find
that
\beq
\frac{E(q^{a+1})^{\lfloor (a+1)/2\rfloor} [q^{a-i};q^{a+1}]_\infty}
            {E(q)}
= 
\frac{E(q^{a+1})^2}{E(q^{(a+1)/2})}
\prod_{\substack{ j=1\\ j\not\in\{a-i,i+1\}}}^{(a-1)/2}
\frac{E(q^{a+1})}{[q^j;q^{a+1}]_\infty} \succeq 0,
\eeq
by \eqn{zq} and \eqn{tcore} with $t=2$ and $q\to q^{(a+1)/2}$.
If $a\equiv 1\pmod{2}$ and $i=(a-1)/2$ we find
\beq
\frac{E(q^{a+1})^{\lfloor (a+1)/2\rfloor} [q^{a-i};q^{a+1}]_\infty}
            {E(q)}
=
\frac{E(q^{(a+1)/2})^{(a+1)/2}}
     {E(q)}
\left(\frac{E(q^{a+1})}
           {E(q^{(a+1)/2})}\right)^{(a-3)/2}
\succeq 0,
\eeq
by \eqn{tcore} since in this case $a\ge3$ and $E(q^2)/E(q)\succeq0$.

\end{proof}

\begin{remark}
Setting $q\to q^{5}$, $a=3$, $k=2$ and $z=q$ in
in \eqn{corgqpi1} we have
\beq
(1-q^8) E(q^5) E(q^{20})^2 \frac{[q^2;q^5]_\infty}{[q;q^5]_\infty^2}
\succeq 0.
\mylabel{eq:crank5a}
\eeq
This leads to an inequality for the crank of partitions mod $5$.
Using \cite[Thm.1]{AG} and \cite[(4.8)]{G88a}
we have
\beq
\sum_{n\ge0} (M(0,5,5n)-M(1,5,5n))q^n 
= E(q^5) \frac{[q^2;q^5]_\infty}{[q;q^5]_\infty^2}.
\mylabel{eq:crank5b}
\eeq
In view of \eqn{crank5a}, we have
\beq
M(0,5,5n) >  M(1,5,5n) \quad \mbox{for $n\ge0$},
\mylabel{eq:crank5c}
\eeq
by checking the result for the first $8$ coefficients.
Here $M(k,t,n)$ is the number of partitions of $n$ with
crank congruent to $k$ mod $t$. This proves \cite[(8.47)]{G88a}
(conjectured in 1988). From \cite[(13)]{E1} we have
\beq
\sum_{n\ge0} (M(2,11,11n+2)-M(1,11,11n+2))q^n
= E(q^{11}) \frac{[q^3;q^{11}]_\infty}{[q,q^4;q^{11}]_\infty}.
\mylabel{eq:crank11a}
\eeq
Setting $q\to q^{11}$, $a=3$, $k=2$ and $z=q$ in
\eqn{corgqpi1} and checking the first $11$ cases we have
\beq
M(2,11,11n+2) >   M(1,11,11n+2) \quad \mbox{for $n\ne3$}.
\mylabel{eq:crank11b}
\eeq
\end{remark}

\begin{prop}
\label{res}
If $\abs{q}<1$ and $z\ne0$ then
\beq
\frac{E(q^2)[z^4;q^2]_\infty}{[z^2;q^2]_\infty [qz^3;q^2]_\infty}
\succeq 0.
\mylabel{eq:res1} 
\eeq
\beq
\frac{E(q^3) (z^2;q^3)_\infty}{ (q^3 z^{-1};q^3)_\infty (z;q)_\infty}
\succeq 0.
\mylabel{eq:res2} 
\eeq
\end{prop}

\begin{proof}
(i) First we find that              
\beq
E(q^3) \,  \frac{[z^2;q^3]_\infty}
                {[z;q]_\infty}        \succeq 0,
\mylabel{eq:res2b}
\eeq
since
\begin{align}
&E(q^3) \,  \frac{[z^2;q^3]_\infty}
                {[z;q]_\infty}  
\mylabel{eq:res2c} \\
& =\left( (1-z^2) \frac{ E(q^3) }
                       {[z;q^3]_\infty}\right)
\left(\frac{(z^2q^3;q^3)_\infty}
           {(zq,zq^2;q^3)_\infty}\right)
\left(\frac{(q^3/z^2;q^3)_\infty}
           {(q/z,q^2/z;q^3)_\infty}\right)
\nonumber,
\end{align}
and each of the three terms on the right side of \eqn{res2c} has nonnegative
coefficients by \eqn{aci} with $m=2$, \eqn{atq} and \eqn{atq} respectively.  
Next we can
use the quintuple product identity to obtain
\beq
\frac{E(q^2)[z^4;q^2]_\infty}
{[z^2;q^2]_\infty [qz^3;q^2]_\infty}
= E(q^6) \frac{[q^2z^6;q^6]_\infty}
              {[qz^3;q^2]_\infty}
+ z^2 E(q^6) \frac{[q^2/z^6;q^6]_\infty}
              {[q/z^3;q^2]_\infty} \succeq 0,
\mylabel{eq:res2d}
\eeq
by \eqn{res2b}.

(ii) We have
\beq
\frac{E(q^3) (z^2;q^3)_\infty}{ (q^3 z^{-1};q^3)_\infty (z;q)_\infty}
=
\left( (1-z^2) \frac{E(q^3)}
                    {[z;q^3]_\infty} \right)
\,
\left(\frac{(z^2q^3;q^3)_\infty}
           {(zq;q^3)_\infty (zq^2;q^3)_\infty}\right) \succeq 0,
\mylabel{eq:res2e}
\eeq
by \eqn{aci}(with $m=2$) and \eqn{atq}.
\end{proof}

\begin{prop}
\label{prop:eta1}
Suppose $n$ is a positive integer.
\begin{enumerate}
\item[(i)]
If $m$ is a positive odd integer then
\beq
\frac{E(q^{mn})^{n(m-1)/2 -m} E(q^m)^{(m+1)/2} E(q^n)}
     {E(q)} \succeq 0.
\mylabel{eq:eta1a}
\eeq
\item[(ii)]
If $m$ is a positive even integer then
\beq
\frac{E(q^{mn})^{(n-2)(m/2-1)} E(q^m)^{m/2-1} E(q^n)E(q^{m/2})}
     {E(q)E(q^{mn/2})} \succeq 0.
\mylabel{eq:eta1b}
\eeq
\end{enumerate}
\end{prop}
\begin{proof}
By \eqn{Cazq} we have
\beq
C_t(z,q) = E(q) E(q^t)^{t-2} \frac{[z^t;q^t]_\infty}
                                  {[z;q]_\infty} \succeq0,
\mylabel{eq:Ctzq}
\eeq
for any positive integer $t$.

(i) Result true for $m=1$ so we suppose $m\ge3$ is an odd integer.
Then
\beq
\prod_{r=1}^{(m-1)/2} [q^r;q^m]_\infty = \frac{E(q)}
                                              {E(q^m)}.
\mylabel{eq:Bev}
\eeq
Hence, by \eqn{Ctzq} we have
\beq
\prod_{r=1}^{(m-1)/2} C_n(q^r,q^m) = 
\frac{E(q^{mn})^{n(m-1)/2 -m} E(q^m)^{(m+1)/2} E(q^n)}
     {E(q)} \succeq 0.
\mylabel{eq:eta1aid}
\eeq

(ii) Result true for $m=2$ so we suppose $m\ge4$ is an even integer.
This time
\beq
\prod_{r=1}^{(m/2-1)} [q^r;q^m]_\infty = \frac{E(q)}
                                              {E(q^{m/2})},
\mylabel{eq:Bodd}
\eeq
and
\beq
\prod_{r=1}^{(m/2-1)} C_n(q^r,q^m) = 
\frac{E(q^{mn})^{(n-2)(m/2-1)} E(q^m)^{m/2-1} E(q^n)E(q^{m/2})}
     {E(q)E(q^{mn/2})} \succeq 0.
\mylabel{eq:eta1bid}
\eeq
\end{proof}
\begin{cor}
\label{cor:eta1}
If $m$ and $n$ are positive integers then
\beq
\frac{E(q^{mn})^{mn-m-n} E(q^m) E(q^n)}
{E(q)} \succeq 0.
\mylabel{eq:eta2}
\eeq
\end{cor}
\begin{proof}
\noindent
We consider two cases.

\noindent
{\it Case 1.} $m$ is odd. By \eqn{eta1a} and \eqn{tcore}
(with $q\to q^m$ and $t=n$) we have
\begin{align}
&\frac{E(q^{mn})^{mn-m-n} E(q^m) E(q^n)}
{E(q)} \mylabel{eta2id1} \\
&=
\left(
\frac{E(q^{mn})^{n(m-1)/2 -m} E(q^m)^{(m+1)/2} E(q^n)}
     {E(q)} 
\right)
\,
\left(
\frac{E(q^{mn})^n}
     {E(q^m)}
\right)^{(m-1)/2}
\succeq 0.
\nonumber
\end{align}

\noindent
{\it Case 2.} $m$ is even. By \eqn{eta1b} (with $n=2$ and $m\to 2n$)
we have
\beq
V_n(q) := 
\frac{E(q^{2n})^{n-2} E(q^2)E(q^{n})}
     {E(q)} \succeq 0,
\mylabel{eq:Vn}
\eeq
where $n$ is any positive integer.
We find that
\begin{align}
&\frac{E(q^{mn})^{mn-m-n} E(q^m) E(q^n)}
{E(q)} \nonumber \\
&=
\left(
\frac{E(q^{mn})^{(n-2)(m/2-1)} E(q^m)^{m/2-1} E(q^n)E(q^{m/2})}
     {E(q)E(q^{mn/2})}
\right)  \mylabel{eq:eta2id2}\\
&\qquad
\left(
\frac{E(q^{mn})^n}
     {E(q^m)}
\right)^{m/2 -1}
\,
V_n(q^{m/2}) \succeq 0,
\nonumber
\end{align}
by \eqn{eta1b}, \eqn{tcore} (with $q\to q^m$ and $t=n$) and
\eqn{Vn}.
\end{proof}

\begin{remark}
We note that in the case when $m$ and $n$ are distinct primes
\eqn{eta2} is a special case of Saito's Conjecture ($N=mn$).
Also, in the case when $m$ is odd there is simple direct
proof. In this case we find that
\beq
0 \preceq \prod_{r=1}^{(m-1)/2} D_n(q^r,q^m)
=
\frac{E(q^{mn})^{mn-m-n} E(q^m) E(q^n)}
     {E(q)},
\mylabel{eq:eta2alt}
\eeq
since each $D_n(q^r,q^m)$ (defined in \eqn{Dazq})
has nonnegative coefficients.
\end{remark}

We make the following
\begin{conj}
\label{conj2}
Suppose $\abs{q}<1$ and $z\ne0$.
\begin{enumerate}
\item[(i)]
If $p\ge1$ then
\beq
\frac{ E(q) } {(z;q)_\infty (qz^{-p};q)_\infty} \succeq 0.
\mylabel{eq:conj2a} 
\eeq
\item[(ii)]
If $a$, $b$, $m$, $n\ge1$ then
\beq
\frac{E(q^{ma + nb})}
{(q^a;q^{ma+nb})_\infty (q^b;q^{ma+nb})_\infty} \succeq 0.
\mylabel{eq:conj2b} 
\eeq
\item[(iii)]
For $n\ge3$                        
\beq
\frac{(z, z^{n-1}q^n;q^n)_\infty}
     {(z;q)_\infty} \succeq 0.
\mylabel{eq:conj2c} 
\eeq
\item[(iv)]
For $n\ge4$ 
\beq
\frac{(z^{n-1}q^n;q^n)_\infty}
     {(zq,zq^2,zq^3;q^n)_\infty} \succeq 0.
\mylabel{eq:conj2d}
\eeq
\item[(v)]
For $n\ge2$                        
\beq
E(q^n) \, \frac{[z^{n-1};q^n]_\infty}
               {[z;q]_\infty} \succeq 0.
\mylabel{eq:conj2e} 
\eeq
\item[(vi)]
For $n>1$, $m>0$, $a=1$, $2$
\beq
\frac{E(q^{nm})}
     {(q^a;q^m)_\infty} \succeq 0.
\mylabel{eq:conj2f}
\eeq
\item[(vii)]
For $m>1$
\beq
E(q^m) \, \frac{[z^2;q^m]_\infty}
            {[z;q^m]_\infty (zq,q/z;q^m)_\infty} \succeq 0.
\mylabel{eq:conj2g}
\eeq
\item[(viii)]
For $n\ge2$
\beq
E(q^n) \frac{[z^{n^2};q^n]_\infty}
            {[z^n;q^n]_\infty}
       \frac{[z^{n+1}q^n;q^n]_\infty}
            {[z^{n+1}q;q]_\infty} \succeq 0.
\mylabel{eq:conj2h}
\eeq
\end{enumerate}
\end{conj}

\begin{remark}
The case $p=1$ of \eqn{conj2a}, the case $m=n=1$ of \eqn{conj2b} 
and the case $a=1$ of \eqn{conj2a} are all special cases
of Proposition \ref{atineq}.
\end{remark}

\begin{remark}
We consider \eqn{conj2c} and let
\beq
P_n(z,q) := 
\frac{(z, z^{n-1}q^n;q^n)_\infty}
     {(z;q)_\infty},
\mylabel{eq:Pnzq}
\eeq
for $n\ge3$. We can show that \eqn{conj2c} holds for
$n=3$, $4$. We have
\beq
P_3(z,q) = \frac{(z^2q^3;q^3)_\infty}
                {(zq,zq^2;q^3)_\infty} \succeq 0,
\mylabel{conj2c3}
\eeq
by \eqn{atq} with $q\to q^3$, $a=zq$ and $t=zq^2$.
Also,
\beq
P_4(z,q) = \frac{(z^3q^4;q^4)_\infty}
                { (zq;q^2)_\infty (zq^2;q^4)_\infty}
= (-zq;q^2)_\infty \frac{(z^3q^4;q^4)_\infty}
                        {(zq^2;q^4)_\infty (z^2q^2;q^4)_\infty}
\succeq 0,
\mylabel{conj2c4}
\eeq
by \eqn{atq} with $q\to q^4$, $a=zq^2$ and $t=z^2q^2$.
\end{remark}

\begin{remark}
When $n\ge4$ it is clear that \eqn{conj2d} implies \eqn{conj2c}.
\end{remark}

\begin{remark}
Finally, we consider \eqn{conj2e}. We observe that
\beq
E(q^n) \, \frac{[z^{n-1};q^n]_\infty}
               {[z;q]_\infty} 
= (1 - z^{n-1})\frac{E(q^n)}
                    {[z;q^n]_\infty}
  \cdot P_n(z,q) P_n(z^{-1},q).
\mylabel{eq:conj2ea}
\eeq
We note that when $n=3$, \eqn{conj2ea} is \eqn{res2c}.
For $n\ge3$, we see that \eqn{conj2c} implies \eqn{conj2e}
by \eqn{aci} with $m=n-1$ and $q\to q^n$. Thus \eqn{conj2e}
holds $n=3$, $4$. It also holds for $n=2$ since
\beq
E(q^2) \, \frac{[z;q^2]_\infty}
               {[z;q]_\infty} 
= \frac{E(q^2)}
       {[zq;q^2]_\infty} \succeq 0,
\mylabel{conj2e2}
\eeq 
by \eqn{zq}.
\end{remark}

\begin{remark}
The case $m=1$ of \eqn{conj2f} is trivial. When
$m=2$ and $a=1$ we have
\beq
\frac{E(q^{2n})}
     {(q;q^2)_\infty}
= \frac{E(q^2)^2}
       {E(q)}
  \cdot
  \frac{E(q^{2n})}
       {E(q^2)}
  \succeq 0,
\mylabel{conj2f2}
\eeq
by \eqn{tcore} with $t=2$.
\end{remark}

\begin{remark}
The case $m=2$ of \eqn{conj2g} is easy:
\beq
E(q^2) \, \frac{[z^2;q^2]_\infty}
            {[z;q^2]_\infty (zq,q/z;q^2)_\infty}               
=
\frac{E(q^2) [z^2;q^2]_\infty}
     {[z;q]_\infty}
=
\frac{E(q^2)}
     {E(q)}
\cdot
(q,-z,-q/z;q)_\infty \succeq 0,
\mylabel{eq:conj2g2}
\eeq
by \eqn{Cazq2}.
\end{remark}

\begin{remark}
The case $n=2$ of \eqn{conj2h} is \eqn{res1}.
It can be shown that \eqn{gqpi} and \eqn{conj2e}
imply \eqn{conj2h}.
\end{remark}

\section{Concluding Remarks}
As noted in the introduction both \eqn{kid} and \eqn{Cazq}
are special cases of Macdonald's identity of type $A$.
It is natural to consider the following questions.

\begin{enumerate}
\item[(i)] Is there a natural analog of $t$-core which extends
\eqn{kid} to other affine root systems?
\item[(ii)] Are there other special cases of Macdonald's
identity for other affine root systems which give nice product
identities analogous to \eqn{Cazq}?
\end{enumerate}

\noindent
\textbf{Acknowledgement}

\noindent
We would like to thank Professor Kyoji Saito for his interest,
comments and for sending copies of his unpublished work \cite{S5},
\cite{SY}. We also would like to thank Robin Chapman, 
Hjalmar Rosengren and Michael Schlosser for helpful discussions.
Finally we thank
Hjalmar Rosengren for the observation
that Theorem \ref{thm1} is a special case of Macdonald's identity
of type $A$ and suggesting the questions above.

\bibliographystyle{amsplain}

\end{document}